\documentclass[dvips,12pt]{article}
\usepackage{latexsym,amsmath,amssymb}
\begin{document}
\numberwithin{equation}{section}
\thispagestyle{empty}
\baselineskip=6.3mm
\begin{center}
{\large{\bf Split general quasi-variational inequality problem}}

\vspace{.3cm}
{ \bf {K.R. Kazmi{\footnote {emails: kr.kazmi.mm@amu.ac.in; krkazmi@gmail.com (K.R. Kazmi)}}}}

\vspace{.3cm}
{\it Department of Mathematics, Aligarh Muslim University, Aligarh-202002, India}\\
\end{center}

\vspace{.2cm}
\parindent=0mm
{\bf Abstract:} In this paper, we introduce a split  general quasi-variational inequality problem which is a natural extension of split  variational inequality problem, quasi-variational and variational inequality problems in   Hilbert spaces. Using projection method,  we propose an  iterative algorithm for the split  general quasi-variational inequality problem and discuss some its special cases. Further, we discuss  the convergence criteria  of these iterative algorithms. The  results  presented in this paper generalize, unify and improve the previously known  many  results for the quasi-variational and variational inequality problems.

\vspace{.25cm}
\parindent=0mm
{\bf Keywords:} Split general quasi-variational inequality problem, Split quasi-variational inequality problem, Split general variational inequality problem, iterative algorithms, convergence analysis.

\vspace{.25cm}

\vspace{.4cm}
\parindent=8mm
\noindent {\bf 1. Introduction}

\vspace{.3cm}
Throughout the paper unless otherwise stated, for each $i\in\{1,2\}$, let $H_i$  be a real Hilbert space  with inner product $\langle \cdot,\cdot \rangle $ and norm $\| \cdot \|$; let $C_i$  be a nonempty, closed and convex subset of $H_i$.

\vspace{.3cm}
The {\it variational inequality problem} (in short, VIP) is to find $x_1\in C_1$  such that 
$$\langle f_1(x_1),y_1-x_1\rangle\geq 0,~~\forall y_1\in C_1, \eqno(1.1)$$
where $f_1:C_1\to H_1$ be a nonlinear mapping. 

\vspace{.3cm} 

Variational inequality theory introduced by Stampacchia [1] and Fichera [2] independently, in early sixties in potential theory and mechanics, respectively, constitutes a significant extension of variational principles. It has been shown that the variational inequality theory provides the natural, descent, unified and efficient framework for a general treatment of a wide class of unrelated linear and nonlinear problem arising in elasticity, economics, transportations, optimization, control theory and engineering sciences, see for instance [3-8]. The development of variational inequality theory can be viewed as the simultaneous pursuit of two different lines of research. On the one hand, it reveals the fundamental facts on the qualitative behavior of solutions to important classes of problems. On the other hands, it enables us to develop highly efficient and powerful numerical methods to solve, for example, obstacle, unilateral, free and moving boundary value problems. In last five decades, considerable interest has been shown in developing various classes of variational inequality problems, both for its own sake and for its applications.

\vspace{.3cm} 
An important generalization of variational inequality problem is quasi-variational inequality problem introduced and studied  by Bensoussan, Goursat and Lions [9] in connection of impulse  control problems. More presicely, for each $i$, let $C_i: H_i \to 2^{H_i}$ be a  nonempty, closed and convex  set valued mapping, where $2^{H_i}$ be the family of all nonempty subsets of $H_i$. The {\it quasi-variational inequality problem} (in short, QVIP) is to find $x_1\in H_1$  such that $x_1\in C_1(x_1)$ and
$$\langle f_1(x_1),y_1-x_1\rangle\geq 0,~~\forall y_1\in C_1(x_1), \eqno(1.2)$$
where $f_1:H_1\to H_1$ be a nonlinear mapping.

\vspace{.30cm}
We observe that if $C_1(x_1)=C_1$ for all $x_1 \in H_1$, then QVIP(1.2) is reduced to VIP(1.1). In many important applications, $C_1(x_1)=m(x_1)+C_1$ for each $x_1 \in H_1$, where $m: H_1 \to H_1$ is a single valued  mapping, see for instance [4,5]. Since then various  generalization of QVIP(1.2)  have been proposed and analyzed,  see for instance [10-14].

\vspace{.30cm} 
Recently, Censor {\it et al.} [15] introduced and studied  the following {\it split  variational inequality problem} (in short, SpVIP): For each $i\in\{1,2\}$, let $f_i:H_i\to H_i$  be a nonlinear mapping  and $A:H_1\to H_2$ be a bounded linear operator with its adjoint operator $A^*$. Then the SpVIP is to:

\vspace{.30cm}
 \noindent Find $x_1^*\in C_1$ such that
$$\langle f_1(x_1^*),x_1-x_1^*\rangle\geq 0,~~\forall x_1 \in C_1,\eqno(1.3a)$$
and such that
$$x_2^*=Ax_1^*\in C_2~~{\rm solves}~\langle f_2(x_2^*),x_2-x_2^*\rangle\geq 0,~~\forall x_2 \in C_2.\eqno(1.3b)$$
SpVIP(1.3a-b) amount to saying: find a solution of variational inequality problem VIP(1.3a), the image of which under a given bounded linear operator  is a solution of VIP(1.3b). It is worth mentioning  that SpVIP is quite general and permits split minimization between two spaces so the image of a minimizer of a given function, under a bounded linear operator, is a minimizer of another function. SpVIP(1.3a-b) is an important generalization of VIP(1.1). It also includes as special case, the split zero problem  and split feasibility problem  which has already been studied and used in practice as a model in intensity-modulated radiation therapy  treatment planning, see [16-18]. For the further related work, we refer to see Moudafi [19],  Byrne {\it et al.} [20], Kazmi and Rizvi [21-24] and Kazmi [25].

 \vspace{.3cm}
 In this paper, we introduce  the following natural generalization of SpVIP(1.3a-b):  For each $i\in \{1,2\}$, let $C_i: H_i \to 2^{H_i}$ be a  nonempty, closed and convex set valued mapping. For each $i\in\{1,2\}$, let $f_i:H_i\to H_i$ and $g_i:H_i\to H_i$ be nonlinear mappings  and $A:H_1\to H_2$ be a bounded linear operator with its adjoint operator $A^*$. Then we consider the problem:

\vspace{.30cm}
 \noindent Find $x_1^*\in H_1$  such that $g_1(x_1^*)\in C_1(x_1^*)$ and
$$\langle f_1(x_1^*),x_1-g_1(x_1^*)\rangle\geq 0,~~\forall x_1 \in C_1(x_1^*),\eqno(1.4a)$$
and such that
$$x_2^*=Ax_1^*\in H_2,~~g_2(x_2^*)\in C_2(x_2^*)~~{\rm solves}~\langle f_2(x_2^*),x_2-g_2(x_2^*)\rangle\geq 0,~~\forall x_2 \in C_2(x_2^*).\eqno(1.4b)$$
We call the problem (1.4a-b), the {\it split  general quasi-variational inequality problem} (in short, SpGQVIP).

\vspace{.30cm}
 Now, we observe some  special cases of SpGQVIP(1.4a-b).

\vspace{.30cm}
 If we set  $g_i=I_i$, where $I_i$ is identity operator on $H_i$, then SpGQVIP(1.4a-b) is reduced to  the following  {\it split  quasi-variational inequality problem} (in short, SpQVIP):

\vspace{.30cm}
 \noindent Find $x_1^*\in H_1$  such that $x_1^*\in C_1(x_1^*)$ and
$$\langle f_1(x_1^*),x_1-x_1^*\rangle\geq 0,~~\forall x_1 \in C_1(x_1^*),\eqno(1.5a)$$
and such that
$$x_2^*=Ax_1^*\in H_2, x_2^*\in C_2(x_2^*)~~{\rm solves}~\langle f_2(x_2^*),x_2-x_2^*\rangle\geq 0,~~\forall x_2 \in C_2(x_2^*),\eqno(1.5b)$$
which appears to be new.

\vspace{.3cm} 

If we set $C_i(x_i)=C_i$ for all $x_i \in H_i$, then SpGQVIP(1.4a-b) is reduced to following {\it split general variational inequality problem} (in short, SpGVIP):
 
\vspace{.30cm}
 \noindent Find $x_1^*\in H_1$  such that $g_1(x_1^*)\in C_1$ and 
$$\langle f_1(x_1^*),x_1-g_1(x_1^*)\rangle\geq 0,~~\forall x_1 \in C_1,\eqno(1.6a)$$
and such that
$$x_2^*=Ax_1^*\in H_2,~~g_2(x_2^*)\in C_2~~{\rm solves}~\langle f_2(x_2^*),x_2-g_2(x_2^*)\rangle\geq 0,~~\forall x_2 \in C_2,\eqno(1.6b)$$
which appears  to be new. 

\vspace{.3cm}
Further, if we set $C_i(x_i)=C_i$ for all $x_i \in H_i$, and $g_i=I_i$, then SpGQVIP(1.4a-b) is reduced to  SpVIP(1.3a-b).

\vspace{.3cm}
Furthermore, if we set $H_2=H_1; C_2(x_2)=C_1(x_1) \forall x_i$;$ f_2=f_1$, and $g_i=I_i$, then SpGQVIP(1.4a-b) is reduced to  QVIP(1.2).
  
\vspace{.3cm}
 Using projection method,  we propose an  iterative algorithm for SpGQVIP(1.4a-b) and discuss some its special cases which are the iterative algorithms for SpQVIP(1.5a-b), SpGVIP(1.6a-b), SpVIP(1.3a-b) and QVIP(1.2). Further, we discuss  the convergence criteria  of these iterative algorithms. The  results  presented in this paper generalize, unify and improve the previously known  many  results for the quasi-variational and variational inequality problems.

\vspace{0.4cm}
\noindent{\bf 2.  Iterative algorithms }

\vspace{0.3cm}
For each $i\in\{1,2\}$, a mapping $P_{C_i}$ is said to be  {\it metric projection} of $H_i$ onto $C_i$ if for every point  $x_i \in H_i$, there exists a unique nearest point in $C_i$ denoted by $P_ {C_i} (x_i)$ such that
$$ \|x_i-P_{C_i}(x_i)\|\leq \|x_i-{\bar{x}}_i\|, ~~ \forall  {\bar{x}}_i \in C_i.$$
 It is well known that $P_{C_i}$ is nonexpansive mapping and satisfies
$$\langle x_i-{\bar{x}}_i ,P_{C_i}(x_i)-P_{C_i}({\bar{x}}_i) \rangle \geq \|P_{C_i}(x_i)-P_{C_i}({\bar{x}}_i)\|^2, ~~\forall x_i,{\bar{x}}_i \in H_i.\eqno(2.1)$$
Moreover, $P_{C_i}(x_i)$ is characterized by:
$$\langle x_i-P_{C_i}(x_i),{\bar{x}}_i-P_{C_i}(x_i) \rangle \leq 0, ~~ \forall  {\bar{x}}_i \in C_i.\eqno(2.2)$$

 \vspace{.3cm}
Further, it is easy to see that the following is true:
$$x_1^* {\rm ~is~a~solution~of~ QVIP(1.2)}\Leftrightarrow x_1^*=P_{C_1(x_1^*)} (x_1^*-\rho_1 f_1(x_1^*)),~~\rho_1>0.$$

Hence, SpGQVIP(1.4a-b) can be reformulated as follows:~~ Find $x_1^* \in H_1 $ with $x_2^*=Ax_1^*$ such that $g_i(x_i^*)\in  C_i(x_i^*)$ and
$$g_i(x_i^*)= P_{C_i(x_i^*)}(g_i(x_i^*)- \rho_i f_i(x_i^*)), \eqno(2.3)$$
for $\rho_i >0$.

\vspace{0.3cm}
Based on above arguments, we propose the following iterative algorithm for approximating a solution to SpGQVIP(1.4a-b).

\vspace{0.3cm}
Let $\{\alpha^n\} \subseteq (0,1)$ be a sequence such that $\sum \limits^{\infty}_{n=1} \alpha^n=+\infty$, and let $\rho_1,~ \rho_2,~ \gamma$ are parameters with positive values.  

\vspace{0.3cm}
\noindent {\bf Iterative algorithm 2.1.} Given $x_1^0\in H_1,$ compute the iterative sequences $\{x_1^n\}$ defined by the iterative schemes:
$$g_1(y^n)= P_{C_1({x_1^n})}(g_1(x_1^n)- \rho_1 f_1(x_1^n)) \eqno(2.4a)$$
$$g_2(z^n)= P_{C_2({Ay^n})}(g_2(Ay^n)- \rho_2 f_2(Ay^n)) \eqno(2.4b)$$
$$x_1^{n+1}=(1-\alpha^n)x_1^n +\alpha^n[y^n+\gamma A^*(z^n-Ay^n)] \eqno(2.4c)$$
for all $n=0,1,2,.....~, \rho_i, \gamma >0.$

\vspace{0.3cm}
	If we set $g_i=I_i$, then Iterative algorithm 2.1 is reduced to following iterative algorithm for SpQVIP(1.5a-b):

\vspace{0.3cm}
\noindent {\bf Iterative algorithm 2.2.} Given $x_1^0\in H_1,$ compute the iterative sequences $\{x_1^n\}$ defined by the iterative schemes:
$$y^n= P_{C_1({x_1^n})}(x_1^n- \rho_1 f_1(x_1^n)) \eqno(2.5a)$$
$$z^n= P_{C_2({Ay^n})}(Ay^n- \rho_2 f_2(Ay^n)) \eqno(2.5b)$$
$$x_1^{n+1}=(1-\alpha^n)x_1^n +\alpha^n[y^n+\gamma A^*(z^n-Ay^n)] \eqno(2.5c)$$
for all $n=0,1,2,.....~, \rho_i, \gamma >0.$

\vspace{0.3cm}
	If we set $C_i(x_i)=C_i$ for all $x_i \in H_i$, then Iterative algorithm 2.1 is reduced to following iterative algorithm for SpGVIP(1.6a-b):

\vspace{0.3cm}
\noindent {\bf Iterative algorithm 2.3.} Given $x_1^0\in H_1,$ compute the iterative sequences $\{x_1^n\}$ defined by the iterative schemes:
$$g_1(y^n)= P_{C_1}(g_1(x_1^n)- \rho_1 f_1(x_1^n)) \eqno(2.6a)$$
$$g_2(z^n)= P_{C_2}(g_2(Ay^n)- \rho_2 f_2(Ay^n)) \eqno(2.6b)$$
$$x_1^{n+1}=(1-\alpha^n)x_1^n +\alpha^n[y^n+\gamma A^*(z^n-Ay^n)] \eqno(2.6c)$$
for all $n=0,1,2,.....~, \rho_i, \gamma >0.$

\vspace{0.3cm}
	If we set $C_i(x_i)=C_i$ for all $x_i \in H_i$, and $g_i=I_i$, then Iterative algorithm 2.1 is reduced to following iterative algorithm for SpVIP(1.3a-b):

\vspace{0.3cm}
\noindent {\bf Iterative algorithm 2.4[25].} Given $x_1^0\in H_1,$ compute the iterative sequences $\{x_1^n\}$ defined by the iterative schemes:
$$y^n= P_{C_1}(x_1^n- \rho_1 f_1(x_1^n))$$
$$z^n= P_{C_2}(Ay^n- \rho_2 f_2(Ay^n)) $$
$$x_1^{n+1}=(1-\alpha^n)x_1^n +\alpha^n[y^n+\gamma A^*(z^n-Ay^n)] $$
for all $n=0,1,2,.....~,\rho_i, \gamma >0.$

\vspace{0.3cm}
If we set $H_2=H_1;~ C_2(x_2)=C_1(x_1)~ \forall x_i$; $f_2=f_1$, and $g_i=I_i$, then Iterative algorithm 2.1 is reduced to following Mann iterative algorithm for QVIP(1.2):

\vspace{0.3cm}
\noindent {\bf Iterative algorithm 2.5.} Given $x_1^0\in H_1,$ compute the iterative sequences $\{x_1^n\}$ defined by the iterative schemes:
$$y^n= P_{C_1}(x_1^n- \rho_1 f_1(x_1^n))$$
$$x_1^{n+1}=(1-\alpha^n)x_1^n +\alpha^n y^n $$
for all $n=0,1,2,.....~,\rho_1>0.$

\vspace{0.3cm}
\noindent {\bf Assumption 2.1.} For all $x_i, y_i, z_i \in H_i$, the operator $P_{C_i({x_i})}$ satisfies the condition:
$$\|P_{C_i({x_i})}(z_i)-P_{C_i({y_i})}(z_i)\| \leq \nu_i\|x_i-y_i\|,$$
for some constant $\nu_i>0$.

\vspace{0.3cm}
\noindent{\bf Definition 2.1.} A nonlinear mapping $f_1:H_1 \to H_1$ is said to be:
\begin{enumerate}
\item[(i)] $\alpha_1$-{\it strongly monotone}, if there exists a constant $\alpha_1 >0$ such that
$$\langle f_1(x)-f_1(\bar{x}), x-\bar{x}\rangle~ \geq \alpha\|x-\bar{x}\|^2, ~~\forall x,\bar{x} \in H_1;$$
\item[(ii)] $\beta_1$-{\it Lipschitz continuous}, if there exists a constant $\beta_1>0$ such that
$$ \|f_1(x)-f_1(\bar{x})\|\leq \beta\|x-\bar{x}\|, ~~\forall x,\bar{x} \in H_1.$$
\end{enumerate}

\vspace{0.4cm}
\noindent {\bf 3. Results}

\vspace{0.3cm}

Now, we study the convergence of the Iterative algorithm 2.1 for SpGQVIP(1.4a-b).

\vspace{0.3cm}
\noindent{\bf Theorem 3.1.} For each $i\in \{1,2\}$, let $C_i: H_i \to 2^{H_i}$ be a  nonempty, closed and convex set valued mapping. Let $f_i:H_i \to H_i$ be $\alpha_i$-strongly monotone with respect to $g_i$  and $\beta_i$-Lipschitz continuous and let $g_i:H_i \to H_i$ be  $\delta_i$-Lipschitz continuous and $(g_i-I_i)$ be $\sigma_i$-strongly monotone, where $I_i$ is the identity operator on $H_i$. Let $A: H_1 \to H_2$ be a bounded linear operator  and $A^*$ be  its adjoint operator. Suppose $x_1^* \in H_1$ is a solution to SpGQVIP(1.4a-b) and Assumption 2.1 hold, then the sequence $\{x_1^n\}$ generated by Iterative algorithm 2.1  converges strongly to $x_1^*$ provided that the constants $ \rho_i $ and $\gamma$ satisfy the following conditions:
$$\left|\rho_1 - \frac{\alpha_1}{\beta_1^2} \right|<\frac{\sqrt{\alpha_1^2-\beta_1^2(\delta_1^2-k_1^2)}}{\beta_1^2}$$
$$\alpha_1 > \beta_1 \sqrt{\delta_1^2-k_1^2};~~ k_1= \left[\frac{\sqrt{2\sigma_1+1}}{1+2\theta_2}-\nu_1\right];~~\delta_1> \left|k_1 \right|;$$
$$0< \theta_2=\frac{1}{\sqrt{2\sigma_2+1}}\left\{\sqrt{\delta_2^2-2\rho_2 \alpha_2 + \rho_2^2 \beta_2^2}+\nu_2\right\};~~\rho_2>0;~~ \gamma \in \left(0, \frac{2}{\|A\|^2} \right)  $$

\vspace{.3cm}
\parindent=8mm
\noindent{\bf Proof.} Since $x_1^* \in H_1$ is a solution to SpGQVIP(1.4a-b), then $x_1^* \in H_1$  be such that $g_i(x_i^*)\in  C_i(x_i^*)$ and
$$g_1(x_1^*)= P_{C_1(x_1^*)}(g_1(x_1^*)- \rho_1 f_1(x_1^*)), \eqno(3.1)$$
$$g_2(Ax_1^*)= P_{C_2(Ax_1^*)}(g_2(Ax_1^*)- \rho_2 f_2(Ax_1^*)), \eqno(3.2)$$
for $\rho_i >0$.

\vspace{.3cm}
From Iterative algorithm 2.1(2.4a), Assumption 2.1 and (3.1), we have
$$\|g_1(y^n)-g_1(x_1^*)\|=\|P_{C_1(x_1^n)}(g_1(x_1^n)- \rho_1 f_1(x_1^n))-P_{C_1(x_1^*)}(g_1(x_1^*)- \rho_1 f_1(x_1^*))\|$$
$$\leq \|P_{C_1(x_1^n)}(g_1(x_1^n)- \rho_1 f_1(x_1^n))-P_{C_1(x_1^n)}(g_1(x_1^*)- \rho_1 f_1(x_1^*))\|$$
$$+\|P_{C_1(x_1^n)}(g_1(x_1^*)- \rho_1 f_1(x_1^*))-P_{C_1(x_1^*)}(g_1(x_1^*)- \rho_1 f_1(x_1^*))\|$$
$$\leq \|g_1(x_1^n)-g_1(x_1^*)- \rho_1 (f_1(x_1^n)- f_1(x_1^*))\|+\nu_1 \|x_1^n-x_1^*\|.$$

Now, using the fact that $f_1$ is $\alpha_1$-strongly monotone with respect to $g_1$ and $\beta_1$-Lipschitz continuous, and $g_1$ is $\delta_1$-Lipschitz continuous, we have
$$\|g_1(x_1^n)-g_1(x_1^*)-\rho_1 (f_1(x_1^n)- f_1(x_1^*))\|^2 \hspace{2.5in}$$
$$ = \|g_1(x_1^n)-g_1(x_1^*)\|^2 -2 \rho_1 \langle f_1(x_1^n)- f_1(x_1^*), g_1(x_1^n)-g_1(x_1^*)\rangle $$
$$+\rho^2 \|f_1(x_1^n)- f_1(x_1^*)\|^2\hspace{2.5in}$$
$$\leq (\delta_1^2-2\rho_1 \alpha_1 + \rho_1^2 \beta_1^2)\|x_1^n-x_1^*\|^2.\hspace{1.8in}$$

As a result, we obtain
$$\|g_1(y^n)-g_1(x_1^*)\| \leq \left\{\sqrt{\delta_1^2-2\rho_1 \alpha_1 + \rho_1^2 \beta_1^2} +\nu_1\right\} \|x_1^n-x_1^*\|. \eqno(3.3)$$

Since $(g_1-I_1)$ is $\sigma_1$-strongly monotone, we have
$$\|y^n-x_1^*\|^2 \leq \|g_1(y^n)-g_1(x_1^*)\|^2 -2\langle (g_1-I_1)y^n-(g_1-I_1)x_1^*, y^n-x_1^* \rangle$$
$$\leq \|g_1(y^n)-g_1(x_1^*)\|^2 -2\sigma_1\|y^n-x_1^*\|^2,\hspace{.5in}$$
which implies
$$\|y^n-x_1^*\| \leq \frac{1}{\sqrt{2\sigma_1+1}}\|g_1(y^n)-g_1(x_1^*)\| \eqno(3.4)$$

From (3.3) and (3.4), we have
$$\|y^n-x_1^*\| \leq \theta_1 \|x_1^n-x_1^*\|, \eqno(3.6)$$
where $\theta_1=\frac{1}{\sqrt{2\sigma_1+1}}\left\{\sqrt{\delta_1^2-2\rho_1 \alpha_1 + \rho_1^2 \beta_1^2} +\nu_1\right\}.$

\vspace{0.3cm}
Similarly, from Iterative algorithm 2.1(2.4b), Assumption 2.1 and (3.2) and using the fact that $f_2$ is $\alpha_2$-strongly monotone with respect to $g_2$ and $\beta_2$-Lipschitz continuous, and $(g_2-I_2)$ is $\sigma_2$-strongly monotone, and $g_2$ is $\delta_2$-Lipschitz continuous, we have
$$\|g_2(z^n)-g_2(Ax_1^*)\| \leq \left\{\sqrt{\delta_2^2-2\rho_2 \alpha_2 + \rho_2^2 \beta_2^2}+\nu_2\right\} \|Ay^n-Ax_1^*\|,\eqno(3.7)$$
and
$$\|z^n-Ax_1^*\| \leq \theta_2 \|Ay^n-Ax_1^*\|, \eqno(3.8)$$
where $\theta_2=\frac{1}{\sqrt{2\sigma_2+1}}\left\{\sqrt{\delta_2-2\rho_2 \alpha_2 + \rho_2^2 \beta_2^2}+\nu_2\right\}.$ 

\vspace{0.3cm}
Next, from Iterative algorithm 2.1(2.4c), we have
$$\|x_1^{n+1}-x_1^*\|\leq  (1-\alpha^n)\|x_1^{n}-x_1^*\|+ \alpha^n[\|y^n-x_1^*-\gamma A^*(Ay^n-Ax_1^*)\|+\gamma\|A^*(z^n-Ax_1^*)\|]\eqno(3.9)$$

Further, using the definition of $A^*$, fact that $A^*$ is a bounded linear operator with $\|A^*\|=\|A\|$, and given condition on $\gamma$, we have
$$\|y^n-x_1^*-\gamma A^*(Ay^n-Ax_1^*)\|^2=\|y^n-x_1^*\|^2-2\gamma \langle y^n-x_1^*,A^*(Ay^n-Ax_1^*)\rangle + \gamma^2\|A^*(Ay^n-Ax_1^*)\|^2$$
$$\leq\|y^n-x_1^*\|^2- \gamma(2- \gamma\|A\|^2) \|Ay^n-Ax_1^*\|^2 $$
$$\leq \|y^n-x_1^*\|^2\hspace{2.1in}\eqno(3.10)$$
and, using (3.8), we have
$$\|A^*(z^n-Ax_1^*)\|\leq \|A\|\|z^n-Ax_1^*\|\hspace{1.5in}$$
$$ \leq \theta_2 \|A\|\|Ay^n-Ax_1^*\|$$
$$\leq \theta_2 \|A\|^2\|y^n-x_1^*\|. \hspace{.2in}\eqno(3.11)$$

Combining (3.10) and (3.11) with inequality (3.9), as a result, we obtain
$$\|x_1^{n+1}-x_1^*\|\leq [1-\alpha^n(1-\theta)] \|x_1^{n}-x_1^*\|,$$
where $\theta= \theta_1(1+\gamma \|A\|^2 \theta_2)$.

\vspace{0.3cm}
Hence, after $n$ iterations, we obtain
 $$\|x_{n+1}-x^*\|\leq \prod\limits^{n}_{j=1}[1-\alpha_j(1-\theta)] \|x_0-x^*\|. \eqno(3.12)$$

It follows from conditions on $\rho_1~\rho_2$ that $\theta \in (0,1)$. Since $\sum \limits^{\infty}_{n=1} \alpha^n=+\infty$ and $\theta \in (0,1)$, it implies in the light of [10] that
$$\lim_{n \to \infty}\prod\limits^{n}_{j=1}[1-\alpha_j(1-\theta)]=0.$$

 Thus it follows from  (3.12) that $\{x_n\}$ converges  strongly  to $x^*$ as $n \to +\infty$. Since $A$ is continuous, it follows from (3.3),(3.6), (3.7) and (3.8) that $y^n \to x_1^*$, $g_1(y^n) \to g_1(x_1^*)$ $Ay^n \to Ax_1^*$, $z^n\to Ax_1^*$ and $g_2(z^n)\to g_2(Ax_1^*)$ as $n \to +\infty$. This completes the proof.\\

If we set $g_i=I_i,$ then Theorem 3.1 reduces to the following result for the convergence of the Iterative algorithm 2.2 for SpQVIP(1.5a-b).

\vspace{0.3cm}
\noindent{\bf Corollary 3.1.} For each $i\in \{1,2\}$, let $C_i: H_i \to 2^{H_i}$ be a  nonempty, closed and convex  set valued mapping. Let $f_i:H_i \to H_i$ be $\alpha_i$-strongly monotone   and $\beta_i$-Lipschitz continuous and let $A: H_1 \to H_2$ be a bounded linear operator  and $A^*$ be  its adjoint operator. Suppose $x_1^* \in H_1$ is a solution to SpQVIP(1.5a-b) and Assumption 2.1 hold, then the sequence $\{x_1^n\}$ generated by Iterative algorithm 2.2  converges strongly to $x_1^*$ provided that the constants $ \rho_i $ and $\gamma$ satisfy the following conditions:
$$\left|\rho_1 - \frac{\alpha_1}{\beta_1^2} \right|<\frac{\sqrt{\alpha_1^2-\beta_1^2(1-k_1^2)}}{\beta_1^2}$$
$$\alpha_1 > \beta_1 \sqrt{1-k_1^2};~~ k_1= \frac{1}{1+2\theta_2}-\nu_1;~~\left|k_1 \right|<1;$$
$$0< \theta_2=\left\{\sqrt{1-2\rho_2 \alpha_2 + \rho_2^2 \beta_2^2}+\nu_2\right\};~~\rho_2>0;~~ \gamma \in \left(0, \frac{2}{\|A\|^2} \right)  $$

If we set $C_i(x_i)=C_i,~ \forall x_i\in H_i$ then Theorem 3.1 reduces to the following result for the convergence of the Iterative algorithm 2.3 for SpGVIP(1.6a-b).

\vspace{0.3cm}
\noindent{\bf Corollary 3.2.} For each $i\in \{1,2\}$, let $C_i$  be a nonempty, closed and convex set in $H_i$. Let $f_i:H_i \to H_i$ be $\alpha_i$-strongly monotone with respect to $g_i$  and $\beta_i$-Lipschitz continuous and let $g_i:H_i \to H_i$ be  $\delta_i$-Lipschitz continuous and $(g_i-I_i)$ be $\sigma_i$-strongly monotone, where $I_i$ is the identity operator on $H_i$. Let $A: H_1 \to H_2$ be a bounded linear operator  and $A^*$ be  its adjoint operator. Suppose $x_1^* \in H_1$ is a solution to SpGVIP(1.6a-b) and Assumption 2.1 hold, then the sequence $\{x_1^n\}$ generated by Iterative algorithm 2.3  converges strongly to $x_1^*$ provided that the constants $ \rho_i $ and $\gamma$ satisfy the following conditions:
$$\left|\rho_1 - \frac{\alpha_1}{\beta_1^2} \right|<\frac{\sqrt{\alpha_1^2-\beta_1^2(\delta_1^2-k_1^2)}}{\beta_1^2}$$
$$\alpha_1 > \beta_1 \sqrt{\delta_1^2-k_1^2};~~ k_1=\frac{\sqrt{2\sigma_1+1}}{1+2\theta_2};~~\delta_1> \left|k_1 \right|;$$
$$0< \theta_2= \sqrt{\frac{\delta_2^2-2\rho_2 \alpha_2 + \rho_2^2 \beta_2^2}{2\sigma_2+1}};~~\rho_2>0;~~ \gamma \in \left(0, \frac{2}{\|A\|^2} \right)  $$

If we set $H_2=H_1;~ C_2(x_2)=C_1(x_1)~ \forall x_i$; $f_2=f_1;~ A=I_1$, and $g_i=I_i$, then Theorem 3.1 reduces to the following result for the convergence of the Iterative algorithm 2.5 for QVIP(1.2).

\vspace{0.3cm}
\noindent{\bf Corollary 3.3.} Let $C_1: H_1 \to 2^{H_1}$ be a  nonempty, closed and convex  set valued mapping. Let $f_1:H_1 \to H_1$ be $\alpha_1$-strongly monotone   and $\beta_1$-Lipschitz continuous. Suppose $x_1^* \in H_1$ is a solution to QVIP(1.2) and Assumption 2.1 hold, then the sequence $\{x_1^n\}$ generated by Iterative algorithm 2.5  converges strongly to $x_1^*$ provided that the constant $ \rho_1 $  satisfies the following conditions:
$$\left|\rho_1 - \frac{\alpha_1}{\beta_1^2} \right|<\frac{\sqrt{\alpha_1^2-\beta_1^2(1-k_1^2)}}{\beta_1^2}$$
$$\alpha_1 > \beta_1 \sqrt{1-k_1^2};~~ k_1= 1-\nu_1;~~\left|k_1 \right|<1.$$

\noindent{\bf Remark 3.1.}~ It is of further research effort to extend the iterative method presented in this paper for solving the split variational inclusions [19] and the split equilibrium problem [22].

\vspace{.5cm}
\noindent{\bf{References}}
\begin{enumerate}

\item [{1.}] Stampacchia, G: Formes bilinearires coercitives sur les ensembles convexes. C.R. Acad. Sci. Paris {\bf{258}}, 4413-4416 (1964) 
\item [{2.}] Fichera, G: Problemi elastostatici con vincoli unilaterali: Il problema di Signorini ambigue condizione al contorno. Attem. Acad. Naz. Lincei. Mem. Cl. Sci. Nat. Sez. Ia {\bf{7}}(8), 91-140 (1963/64)
\item [{3.}] Bensoussan, A, Lions, JL: Applications of Variational Inequalities to Stochastic Control. North-Holland, Amsterdam, 1982
\item [{4.}] Bensoussan, A, Lions, JL: Impulse Control and Quasivariational Inequalities. Gauthiers Villers, Paris, 1984.

\item [{5.}] Baiocchi, C, Capelo, A: Variational and Quasi-variational Inequalities. Wiley, New York, 1984

\item [{6.}] Crank, J: Free and Moving Boundary Problems. Clarendon Press, Oxford, 1984

\item [{7.}] Glowinski, R:  Numerical Methods for Nonlinear Variational Problems. Springer, Berlin, 1984
\item [{8.}] Kikuchi, N, Oden, JT: Contact Problems in Elasticity. SIAM, Philadelphia, 1998
\item [{9.}] Bensoussan, A, Goursat,M, Lions, JL: Contr\^{o}le impulsinnel et inequations quasivariationnelles stationeries. C.R. Acad. Sci. {\bf 276}, 1279-1284 (1973)
\item [{10.}] Kazmi, KR: On a class of quasivariational inequalities. New Zealand J.  Math. {\bf 24}, 17-23  (1995) 
\item [{11.}] Kazmi, KR: Mann and Ishikawa type perturbed iterative algorithms for generalized quasi-variational inclusions. J. Math. Anal. Appl. {\bf 209}, 572-584 (1997)
\item [{12.}] Kazmi, KR, Bhat, MI, Khan, FA: A class of multi-valued quasi-variational inequalities. J. Nonlinear  Convex Anal. {\bf 6}(3), 487-495 (2005) 
\item [{13.}]Kazmi, KR: Iterative algorithm for a class of generalized quasi-variational inclusions with fuzzy mappings in Banach spaces. J. Comput. Appl. Math. {\bf  188}(1), 1-11  (2006) 

\item [{14.}]	Kazmi, KR, Khan, FA, Shahzad, M:  Existence and iterative approximation of a unique solution of a system of general quasi-variational inequality problems, Thai J. Math. {\bf  8}(2), 405-417 (2010)

\item [{15.}]   Censor, Y,  Gibali, A, Reich, S: Algorithms for the  split variational inequality problem. Numerical Algorithms {\bf 59}, 301-323 (2012)
\item [{16.}]   Censor, Y, Bortfeld, T, Martin, B, Trofimov, A: A unified approach for inversion problems in intensity modulated radiation therapy. Physics in Medicine and Biology  {\bf 51},  2353-2365 (2006)
\item [{17.}]  Censor, Y,  Elfving, T: A multiprojection  algorithm using Bergman projections in product space. Numerical Algorithms {\bf 8}, 221-239 (1994)
\item [{18}] Combettes, PL: The convex feasibility problem in image recovery. Adv. Imaging Electron Phys. {\bf 95}, 155-270 (1996)

\item[{19.}]  Moudafi, A: Split monotone variational inclusions. J. Optim. Theory Appl. {\bf 150}, 275-283 (2011)
\item[{20.}] Byrne, C, Censor, Y, Gibali, A, Reich, S:  Weak and strong convergence of algorithms for the  split common null point problem. J. Nonlinear Convex Anal. {\bf 13}, 759-775 (2012)
\item[{21.}]  Kazmi, KR, Rizvi, SH: Iterative approximation of a common solution of a split equilibrium problem, a variational inequality problem and a fixed point problem.  J. Egyptian Math. Soc. {\bf 21}, 44-51 (2013) 
\item[{22.}] Kazmi, KR, Rizvi, SH: Iterative approximation of a common solution of a split generalized equilibrium problem and a fixed point problem for nonexpansive semigroup, Mathematical Sciences {\bf 7}, Art. 1 (2013) (doi 10.1186/2251-7456-7-1)
\item[{23.}] Kazmi, KR, Rizvi, SH: An iterative method for split variational inclusion problem and fixed point problem for a nonexpansive mapping. In Press, Optimization Letters (2013)(doi 10.1007/s11590-013-0629-2)
\item[{24.}] Kazmi, KR, Rizvi, SH: Implicit iterative method for approximating a common solution of split equilibrium problem and  fixed point problem for a nonexpansive semigroup. In Press,  Arab J. Math. Sci. (2013)(doi 10.1016/j.ajmsc.2013.04.002)

\item[{25.}]  Kazmi, KR: Split nonconvex variational inequality problem. Mathematical Sciences  {\bf 7}, Art. 36 (2013) (doi: 10.1186/10.1186/2251-7456-7-36)

\end{enumerate}

\end{document}